%% file: bs.tex
\documentclass[12pt]{article}
\usepackage{amsmath}
\usepackage{amssymb}

\usepackage{epsfig}

\newtheorem{theo}{Theorem}[section]
\newtheorem{lmm}{Lemma}[section]

\allowdisplaybreaks

\numberwithin{equation}{section}

\def\Z{{\mathbb Z}}

\def\Pr{{\mathbf P}}
\def\1{{\mathbf 1}}

\def\MM{{\mathcal M}}
\def\ZZ{{\mathcal Z}}
\def\RR{{\mathcal R}}
\def\X{{\mathfrak X}}
\def\T{{\mathcal T}}
\def\VV{{\mathcal V}}
\def\D{{\mathfrak D}}
\def\tD{{\tilde{\mathfrak D}}}
\def\PP{{\mathcal P}}
\def\tPP{{\tilde{\mathcal P}}}
\def\L{{\mathcal L}}
\def\H{{\mathcal H}}

\def\hT{{\hat{\mathcal T}}}

\def\kk{{\kappa_N}}

\def\th{{\tilde h}}
\def\tV{{\tilde V}}
\def\tE{{\tilde E}}
\def\hh{{\hat h}}
\def\ha{{\hat a}}
\def\tQ{{\tilde Q}}

\newcommand{\ord}{\mathop{\sf ord}}

\newcommand{\lc}{\lceil}
\newcommand{\rc}{\rceil}

\newcommand{\lf}{\lfloor}
\newcommand{\rf}{\rfloor}

\newcommand{\sq}{N^{1/2}}

\def\eps{\varepsilon}
\let\phi=\varphi
\def\qed{\hfill\rule{.2cm}{.2cm}}

\begin{document}

\title{Random Bulgarian solitaire}

\author{Serguei~Popov\thanks{Partially supported by CNPq (302981/02--0)}}

\maketitle

{\footnotesize
\noindent
Departamento de Estat\'\i stica, Instituto de Matem\'atica e
Estat\'\i stica, Universidade de S\~ao Paulo, rua do Mat\~ao 1010,
CEP 05508--090, S\~ao Paulo SP, Brasil.\\
E-mail: popov@ime.usp.br

}

\begin{abstract}
We consider a stochastic variant of the game of Bulgarian
solitaire~\cite{G}. For the stationary measure of the random
Bulgarian solitaire, we prove that most of its mass is
concentrated on (roughly) triangular configurations of
certain type.
\\[0.3cm]
{\bf Keywords:} shape theorem, triangular configuration,
Markov chain, stationary measure
\end{abstract}

\section{Introduction and results}
\label{s_introres} Consider the following (random) game: a deck
of~$N$ cards is divided into several piles.
Then, for each pile, we leave it intact with
probability~$1-p$ and remove
one card from there with probability~$p$
($p\in [0,1]$ is a given parameter),
independently of other piles. The cards that were removed are
collected to form a new pile. The order of piles is not important
and the piles of size zero are ignored. The case $p=0$ is trivial
(nothing moves) and will not be considered. When $p=1$, this is the
game of {\it Bulgarian solitaire}, made known by Martin Gardner in~\cite{G},
and studied in \cite{AD,Be,GH,E,I} (cf.\ also~\cite{CH,Y} for
some variations of that game).
The ``truly random'' model with parameter $0<p<1$ is a discrete-time
irreducible and aperiodic Markov chain on the space of all unordered
partitions of~$N$; for obvious reasons, it will be referred to as
the {\it random Bulgarian solitaire}.

If the number of cards~$N$
is a triangular number, i.e., $N=1+2+\cdots+k$ for some~$k$,
a remarkable fact is that, starting from any initial configuration,
after a finite number of moves the (deterministic) Bulgarian solitaire
will reach the stable configuration formed by piles of sizes $k,k-1,\ldots,1$.
The above result was proved in~\cite{VGRT} (see the solution to
Problem 6.10)
and in~\cite{Br} independently, and later it was discovered
that the maximal number of moves necessary to enter the stable
configuration is $k^2-k$, and that that bound is sharp (see~\cite{E,I}).
If~$N$ is not a triangular number, then such a stable configuration does not
exist. However, it is possible to prove that after at most $O(k^2)=O(N)$
moves the game will enter into a cycle. Moreover, all the configurations
of the cycle are ``almost triangular'' in the following sense. Let
$k=\max\{n: 1+2+\cdots+n\leq N\}$; then all the configurations in that
cycle can be constructed from the configuration $(k,k-1,\ldots,1)$ by
adding at most one card to each pile,
and maybe adding one more pile of size~$1$, see \cite{AD,Be,GH,E} for
exact formulations and more details.

Thus, we see that Bulgarian solitaire ``likes'' triangular configurations,
and so we may expect some kind of similar behaviour from the random
Bulgarian solitaire. There is no possibility, however, to obtain exact results
of the form of those of the previous paragraph,
since random Bulgarian solitaire
is a finite irreducible Markov chain, so it visits all its states
infinitely many times a.s. Instead, we aim at the
results of the following kind:
the stationary measure of the set of configurations which are in some sense
close to the (rescaled) triangular configuration is
close to~$1$.
This can be regarded as a ``shape theorem'' result even though it is
substantially different from most of the shape results appearing
in the literature. (In most cases some time-dependent random
set is constructed, and then, when rescaled by time, it converges to
some, usually nonrandom, shape. See e.g.~\cite{AMP,DL,R,Z} for results
of this kind.) The results we are aiming at resemble rather those
of~\cite{CEP,CLP}.

Also, let us remark here
that the question of how fast the deterministic Bulgarian solitaire
{\it approximates\/} the triangle has not been yet studied in the
literature. To motivate this question, take $N=1+2+\cdots+k$, and suppose
that the initial configuration is $(k-1,k-1,k-2,k-3,\ldots,3,2,1,1)$, i.e.,
the exact triangular configuration is modified by removing one card from the
biggest pile and forming one more pile of size~$1$ with that card.
Then, macroscopically this configuration is already quite triangular;
however, if we are aiming to reach $(k,k-1,\ldots,1)$, this is the worst
possible initial configuration (the number of moves needed is exactly $k^2-k$,
cf.~\cite{E})!
Here we prove that, whenever the initial configuration is ``reasonable''
(i.e., the number of piles is~$O(\sq )$ and the number of cards in the biggest
pile is also~$O(\sq )$), we need only~$O(\sq )$ moves
of deterministic Bulgarian solitaire to make the
($\sq $-rescaled) configuration {\it close\/} to the triangle.
While such a result by itself may not be of great interest,
the method of its proof will be an important tool in the course of
the proof of the results about random Bulgarian solitaire.

Now, we introduce some notations and give the formal definition of the
process. If $\ell(S)$ is the number of piles in the configuration~$S$,
we write $S=(k_1,\ldots,k_{\ell(S)})$, where $k_1\geq\ldots\geq k_{\ell(S)}$.
We denote also by $R(S):=k_1$ the size of the biggest pile and
by $|S|:=k_1+\cdots+k_{\ell(S)}$ the number of cards in the configuration.
Let $\ord(n_1,\ldots,n_m)$ be the operation of
putting $n_1,\ldots,n_m$ in the decreasing order and discarding zeros. Now, let
$\xi_1,\xi_2,\xi_3,\ldots$ be a sequence of i.i.d.\ random variables such that
$\Pr[\xi_1=1]=1-\Pr[\xi_1=0]=p$. Then the operator $Q_p$ which transforms
the configuration $S=(k_1,\ldots,k_{\ell(S)})$
in the game of random Bulgarian solitaire with parameter~$p$
is defined by
\[
Q_p S = \ord(k_1-\xi_1,\ldots,k_{\ell(S)}-\xi_{\ell(S)},
                                   \xi_1+\cdots+\xi_{\ell(S)}).
\]
Denote also by~$Q_p^{(n)}S$ the result of~$n$ independent applications
of~$Q_p$ to~$S$; clearly, the process is conservative in the sense
that $|Q_p^{(n)}S|=|S|$ for all~$n$. Suppose that
$|S_0|=N$. As remarked above, for $0<p<1$ the stochastic process
$S_0,Q_pS_0,Q_p^{(2)}S_0,\ldots$ is an irreducible aperiodic Markov chain with
finite state space $\X_N:=\{S:|S|=N\}$. We denote by~$\pi_{p,N}(\cdot)$
its stationary measure.

To formulate our results, we need also to find a way to define sets
of configurations that are close to a specific triangular configuration.
To this end, for two configurations~$S_1=(n_1,\ldots,n_{\ell(S_1)})$,
$S_2=(m_1,\ldots,m_{\ell(S_2)})$ define the distance~$\rho(S_1,S_2)$ by
\[
\rho(S_1,S_2) = \max_{j\geq 1}|n_j-m_j|,
\]
with the convention $n_j=0$ for all $j>\ell(S_1)$ and
$m_j=0$ for all $j>\ell(S_2)$. Next, we define the triangular
configuration $\T(p,N)=(n_1,\ldots,n_{m_0})$ by
$n_j=\lc (2Np)^{1/2} - pj\rc$, $m_0=\ell(\T(p,N))=\max\{j:
    \lc (2Np)^{1/2} - pj\rc \geq 1 \}$. When~$p$ is fixed and~$N\to\infty$,
    we can write $R(\T(p,N))=(2Np)^{1/2}+O(1)$,
    $\ell(\T(p,N))=(2N/p)^{1/2}+O(1)$. Finally, for~$\eps>0$
    (which may depend on~$N$) define the set $\T(\eps,p,N)$ of
``roughly triangular'' configurations by
\[
\T(\eps,p,N) = \{S: |S|=N, \rho(S,\T(p,N))\leq \eps \sq \}.
\]

Let $k=k(N)=\max\{n: 1+2+\cdots+n\leq N\}$, and define
the configuration $\T_0^N:=(k+m_1,k-1+m_2,\ldots,1+m_k)$,
where
\[
m_i = \left\{
         \begin{array}{ll}
          1, & \mbox{ if } i\leq N - \frac{k(k+1)}{2},\\
          0, & \mbox{ otherwise.}
         \end{array}
         \right.
\]
Note that $|\T_0^N|=N$ (for example, for $N=11$ we have $\T_0^N=(5,3,2,1)$).
For the particular case~$p=1$ we say that $\T(\eps,1,N)$ is
nondegenerate if it contains the configuration
$\T_0^N$, as well as all the configurations~$S$
with $\rho(\T_0^N,S)=1$.
It is easy to see that for any fixed~$\eps>0$ there exists $N_0=N_0(\eps)$
such that $\T(\eps,1,N)$ is nondegenerate for all $N\geq N_0$, and
the same is true when e.g.\ $\eps\sim N^{-\alpha}$, $\alpha<1/2$.

Now we are ready to formulate the main results of this paper.
First, we state the result about the time to approximate the triangular
configuration for the deterministic Bulgarian solitaire (i.e., with $p=1$).
\begin{theo}
\label{det_BS}
Take $\eps>0$ and suppose that~$N$ is large enough
to guarantee that $\T(\eps,1,N)$ is nondegenerate.
Suppose that the initial configuration~$S_0$ with $|S_0|=N$ has the following
properties: $\ell(S_0)\leq \gamma_1 \sq$ and $R(S_0)\leq \gamma_2\sq$ for some
$\gamma_1,\gamma_2>0$. Then there exists
$v_0=v_0(\eps,\gamma_1,\gamma_2)$ such that we have
\begin{equation}
\label{eq_t1}
 Q_1^{(n)}S_0 \in \T(\eps,1,N)
\end{equation}
for all $n\geq v_0\sq$.
\end{theo}
In words, this result means that if the initial configuration is ``reasonable'',
then the number of moves required to approximate the triangle is $O(\sq)$.

Now, we turn our attention to random Bulgarian solitaire:
\begin{theo}
\label{ran_BS}
Suppose that $0<p<1$. Then
for any $a<1/144$ there exist positive constants
 $v_1=v_1(a,p)$ and $\delta=\delta(a,p)$
such that for all~$N$
\begin{equation}
\label{eq_t2}
 \pi_{p,N}(\T(N^{-a},p,N)) \geq 1-\exp(-v_1 N^{\delta}).
\end{equation}
\end{theo}

In Section~\ref{s_fr} there are some more comments and open problems
related to the Bulgarian solitaire (both deterministic and random).
Also, the reader may find it interesting to look at JAVA simulation of the random
Bulgarian solitaire (with $p=1/2$) on the internet page of Kyle Petersen at\\
{\tt http://people.brandeis.edu/\~{}tkpeters/reach/stuff/reach}

\section{Proofs}
\label{s_proofs}
This section is organized in the following way. In Section~\ref{s_Etienne} we
introduce the notion of {\it Etienne diagram}, which is just
another way to represent the configurations of the game.
Then, we show how the moves of Bulgarian solitaire are performed on this
diagram and discuss its other properties. In Section~\ref{s_proof_det}
we prove Theorem~\ref{det_BS}, and in Section~\ref{s_proof_rand}
we prove Theorem~\ref{ran_BS} (in Section~\ref{s_proof_rand} some
results and technique from Sections~\ref{s_Etienne} and~\ref{s_proof_det}
are used, most notably the inequality~(\ref{iteration})).

\subsection{Representation via Etienne diagram and its properties}
\label{s_Etienne}
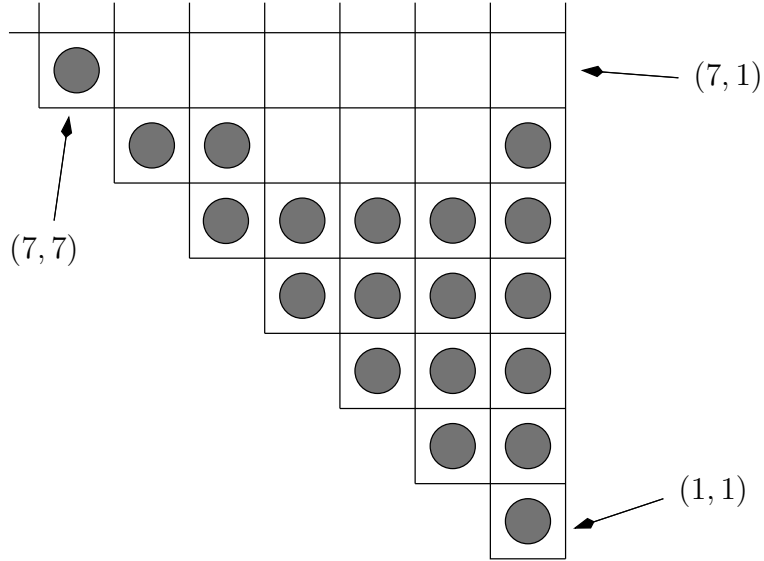
\begin{figure}
\begin{center}
\input{p1.pstex_t}
\caption{The Etienne diagram of $S=(7,5,3,2,1,1)$. We have $\RR_{i,j}(S)=1$
for all $(i,j)$ with $i\leq 5$ and for $(i,j)=(6,1), (6,5), (6,6), (7,7)$}
\label{p1}
\end{center}
\end{figure}

Before starting the proofs, we need to describe another representation
of a particular state of (deterministic)
Bulgarian solitaire, which we call an Etienne diagram
(cf.~\cite{E}). In this approach the cards are identified with particles
living in the cells of the set $\ZZ=\{(i,j)\in\Z^2: i\geq 1, 1\leq j \leq i\}$,
with at most one particle per cell. We write $\RR_{i,j}=1$ when
the cell $(i,j)$ is occupied and $\RR_{i,j}=0$ when the cell $(i,j)$ is
empty. Clearly, $\ZZ$ is a half-quadrant of~$\Z^2$, but we would like
to visualize~$\ZZ$ in a little bit unconventional way (see Figure~\ref{p1}):
the cell~$(1,1)$ lies in the base and supports the column
$\{(i,1), i=1,2,3,\ldots\}$,
while the diagonal $\{(i,i), i=1,2,3,\ldots\}$ goes in the NW direction
(so the rows of~$\ZZ$ are enumerated from right to left; notice that at this
point we deviate from~\cite{E},
where the rows were enumerated from left to right).
Now, a configuration~$S=(k_1,\ldots,k_{\ell(S)})$ is represented
as follows (as on Figure~\ref{p1}): we put $\RR_{i,j}=\RR_{i,j}(S)=1$ for
\[
  (i,j) \in \bigcup_{n=1}^{\ell(S)}\bigcup_{m=1}^{k_n}\{(n+m-1,m)\},
\]
and $\RR_{i,j}=\RR_{i,j}(S)=0$ for all other pairs $(i,j)$. From
the fact that $k_1\geq\ldots\geq k_{\ell(S)}$ we immediately deduce
that for any~$S$
\begin{equation}
\label{pust_tri}
\mbox{if $\RR_{i,j}=0$ then $\RR_{n,m}=0$ for all
$n\geq i, j\leq m\leq j+n-i$},
\end{equation}
and
\begin{equation}
\label{poln_tri}
\mbox{if $\RR_{i,j}=1$ then $\RR_{n,m}=1$ for
all $n\leq i, \max\{1,j-i+n\}\leq m\leq j$}.
\end{equation}

One of the advantages of the representation via Etienne diagram is that
it makes it more clear how the process approaches the triangular
configuration. To see what we mean, first note that the move of the Bulgarian
solitaire consists in applying the following two substeps to the
corresponding Etienne diagram (see Figure~\ref{p2}):
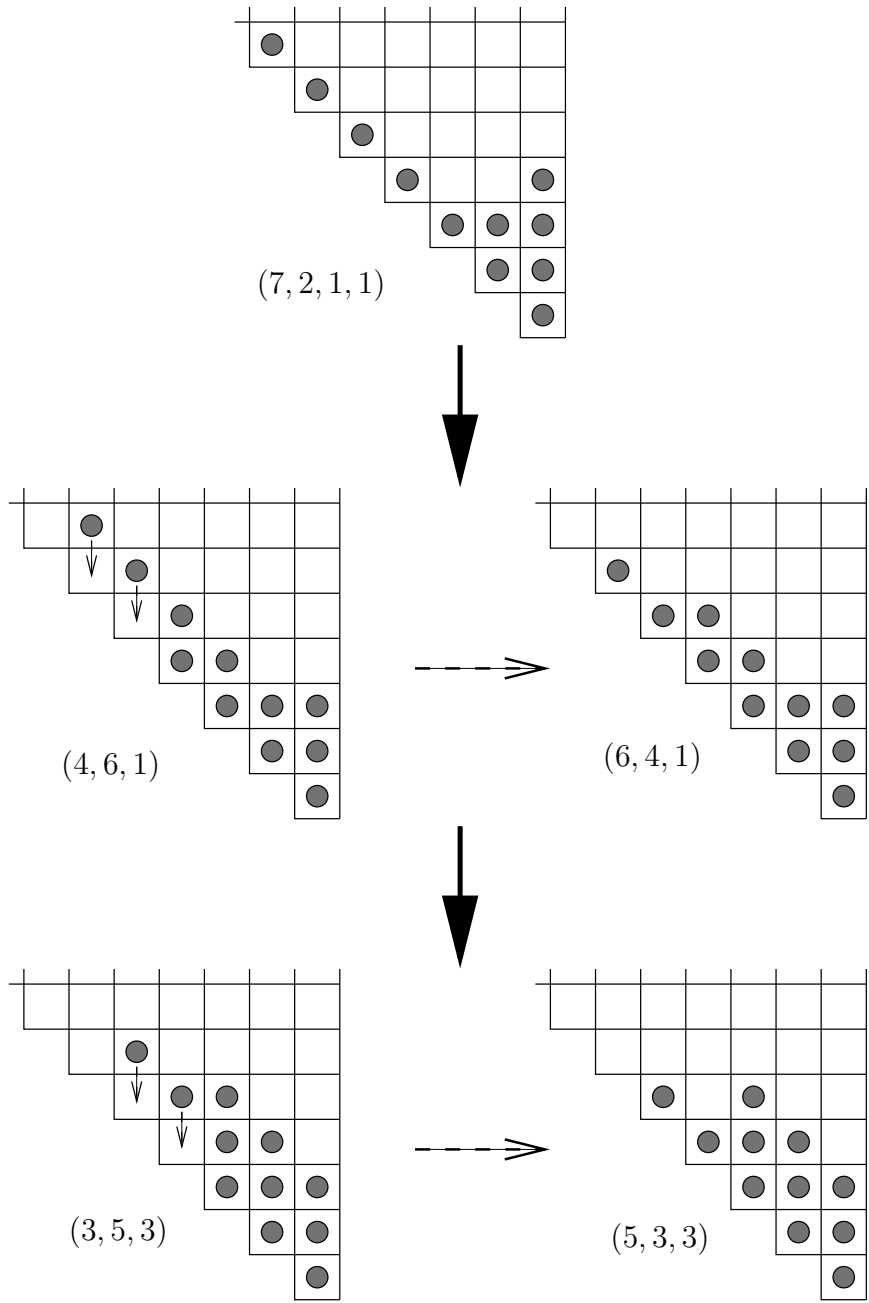
\begin{figure}
\begin{center}
\input{p2.pstex_t}
\caption{Moves of Bulgarian solitaire on the Etienne diagram}
\label{p2}
\end{center}
\end{figure}
\begin{itemize}
\item apply the cyclic shift (from left to right)
   to each row of the diagram;
\item if after the shift there is a particle that is placed above an empty cell,
   then the particle falls there; this procedure is repeated until no further
fall is possible.
\end{itemize}
Speaking formally, let $S'=Q_1S$. Then the Etienne diagram
of~$S'$ is constructed
using the following procedure:
\begin{itemize}
\item[{\bf (I):}] For all~$i$ put $b^0_{i,j}=\RR_{i,j+1}(S)$ for $j<i$
   and $b^0_{i,i}=\RR_{i,1}(S)$.
\item[{\bf (II):}] Suppose that for the array~$b^0$
we can find $(i_0,j_0)$ such that
$b^0_{i_0,j_0}=0$, $b^0_{i_0+1,j_0}=1$. Then construct the new array~$b^1$
by $b^1_{i_0,j_0}=1$, $b^1_{i_0+1,j_0}=0$, and $b^1_{i,j}=b^0_{i,j}$
for $(i,j)\neq(i_0,j_0),(i_0+1,j_0)$.
\item[{\bf (III):}] Repeat the previous procedure with~$b^1$
instead of~$b^0$, and so on.
At some moment we will obtain an array~$b^{\tilde m}$ for which we cannot
find $(i_{\tilde m},j_{\tilde m})$ such that
$b^{\tilde m}_{i_{\tilde m},j_{\tilde m}}=0$,
$b^{\tilde m}_{i_{\tilde m}+1,j_{\tilde m}}=1$. Then for all $(i,j)\in\ZZ$ put
$\RR_{i,j}(S')=b^{\tilde m}_{i,j}$.
\end{itemize}

Now, suppose that $|S|=1+2+\cdots+k$. On the Etienne diagram the triangular
configuration corresponds to the configuration $(\RR_{i,j}=\1\{i\leq k\})$.
Note also that if the first~$m$ rows of the diagram are occupied, then
they will remain occupied during all the subsequent evolution. This shows that
the falls of particles ``help'' to reach the stable configuration
(more and more rows become all occupied). Moreover, in many concrete situations
it is possible to know how many moves are needed to fill out some
region which was originally empty. Arguments of this kind will be heavily used
in the course of the proof of our results.

Consider the Etienne diagram of a configuration~$S$.
Since the system is conservative, there is a natural
correspondence between particles (holes) in that diagram
and particles (holes) in the diagram of the configuration~$Q_1S$.
This shows that for each particle (hole) on the original diagram
we can define its trajectory, i.e., we know its position after~$n$
moves of the game. Let $J_{i,j}(n)$ be the second coordinate of
the particle (hole) from~$(i,j)$ after~$n$ moves, and let~$\MM_{i,j}(n)$
be the number of falls (movements upwards) that the particle (hole) from
$(i,j)$ was subjected to during~$n$ moves.
That means that, if $\RR_{i,j}(S)=1$, then $(i-\MM_{i,j}(n),J_{i,j}(n))$
are the coordinates of the particle from~$(i,j)$ after~$n$ moves,
while if $\RR_{i,j}(S)=0$, then $(i+\MM_{i,j}(n),J_{i,j}(n))$
are the coordinates of the hole from~$(i,j)$ after~$n$ moves.
It seems to be very
difficult to calculate exactly $J_{i,j}(n)$ and $\MM_{i,j}(n)$
(except in trivial situations, when, e.g., $\RR_{i,j}(S)=1$ and
$\RR_{i',j'}(S)=1$ for all~$i'<i$). However, we can establish some
relation between these quantities by defining first
\[
 {\hat J}_{i,j}(n) = \left\{
      \begin{array}{ll}
       j-n + i\Big\lf
                 \displaystyle\frac{n}{i} \Big\rf \vphantom{\int\limits_{N_N}},
                  & \mbox{if }
                               n-i\Big\lf \displaystyle\frac{n}{i} \Big\rf<j,\\
       j-n + i\Big(\Big\lf
                \displaystyle\frac{n}{i} \Big\rf+1\Big), & \mbox{if }
                             n-i\Big\lf \displaystyle\frac{n}{i} \Big\rf\geq j.
      \end{array}
     \right.
\]
In words, $(i,{\hat J}_{i,j}(n))$ would be the position of particle (hole)
 from~$(i,j)$ at time~$n$ if we know that $\MM_{i,j}(n)=0$
 (the quantities $J,\MM,\hat J$ depend also on~$S$, but we do not
 indicate that in our notations).
\begin{lmm}
\label{JM}
If $\RR_{i,j}(S)=0$ and~$n$ is such that
$i-{\hat J}_{i,j}(n)>\lf n/i\rf\MM_{i,j}(n)$, or
$\RR_{i,j}(S)=1$ and~$n$ is such that
${\hat J}_{i,j}(n)>\lf n/i\rf\MM_{i,j}(n)$,
then
\[
 |J_{i,j}(n)-{\hat J}_{i,j}(n)| \leq \Big\lf \frac{n}{i}\Big\rf\MM_{i,j}(n).
\]
\end{lmm}

\noindent
{\it Proof.}
Suppose for example that $\RR_{i,j}(S)=0$.
Denote $j'=j-\MM_{i,j}(n)$. Since~$n$ is such that
$i-{\hat J}_{i,j}(n)>\lf n/i\rf\MM_{i,j}(n)$,
we have that ${\hat J}_{i,j}(n)\leq {\hat J}_{i,j'}(n)$.
The lemma then follows from the fact that $J_{i,j}(n)$ should be
somewhere in between ${\hat J}_{i,j}(n)$ and ${\hat J}_{i,j'}(n)$.
 The other case is treated analogously.
\qed

Next, we define some quantities which concern the geometric structure
of the representation via Etienne diagram,
and prove some relations between them.
For~$N\geq 1$ define
\[
   \theta_N = \max\Big\{k: \frac{k(k+1)}{2}\leq N\Big\};
\]
when~$N\to \infty$, we have $\theta_N=(2N)^{1/2}+O(1)$.
Using the Etienne representation of a configuration~$S$, define
\begin{eqnarray}
E_-(S) &=& \sum_{i\leq \theta_{|S|}}\sum_{j:\RR_{i,j}(S)=0}
          \Big(\theta_{|S|}-i+\frac{1}{2}\Big),
\label{def_E-}\\
E_+(S) &=& \sum_{i > \theta_{|S|}}\sum_{j:\RR_{i,j}(S)=1}
               \Big(i-\theta_{|S|}-\frac{1}{2}\Big),
\label{def_E+}
\end{eqnarray}
and put $E(S)=E_-(S)+E_+(S)$.
The quantity~$E(S)$ can be thought of as the ``energy'' of the configuration:
the bigger is~$E(S)$, the ``more distant'' (not necessarily in the sense of
the distance~$\rho$) is~$S$ from $\T_0^N$.
Denote also
$G^N_{\alpha,\beta} = \{S: |S|=N, \ell(S)\leq \alpha\sq, R(S)\leq\beta\sq\}$.
The next lemma establishes some elementary properties of the energy~$E(S)$.
\begin{lmm}
\label{prop_E}
\begin{itemize}
\item[(i)] There exists a constant $\gamma=\gamma(\alpha,\beta)$ such that
for all~$N$ and all $S\in G^N_{\alpha,\beta}$ we have
$E(S)\leq \gamma N^{3/2}$.
\item[(ii)] For all~$S$ it holds that $E(Q_1S)\leq E(S)$. Moreover,
$E(S)-E(Q_1S)$ is equal to the number of falls of particles during
the second substep of the move of the Bulgarian solitaire represented
by the Etienne diagram (i.e., it is equal to~$\tilde m$ in {\bf (III)}).
\end{itemize}
\end{lmm}

\noindent
{\it Proof.}
Define $\H(S)=\max\{i: \mbox{ there exists~$j$ such that }\RR_{i,j}(S)=1\}$.
   From (\ref{poln_tri}) one easily gets that there exists
$\gamma'=\gamma'(\alpha,\beta)$ such that
for all $S\in G^N_{\alpha,\beta}$ we have
$\H(S)\leq \gamma'\sq$. The proof of~(i) then reduces to an
elementary computation (roughly speaking, to compute~$E(S)$ we
have at worst~$O(N)$ terms, each of order~$O(\sq)$).

As for the proof of~(ii), note first that the
operation of cyclic shift does not
change the quantities defined in~(\ref{def_E-})--(\ref{def_E+}). Then, it is
straightforward to see that each particle's fall decreases~$E$ by one unit,
which concludes the proof of the lemma.
\qed
\medskip

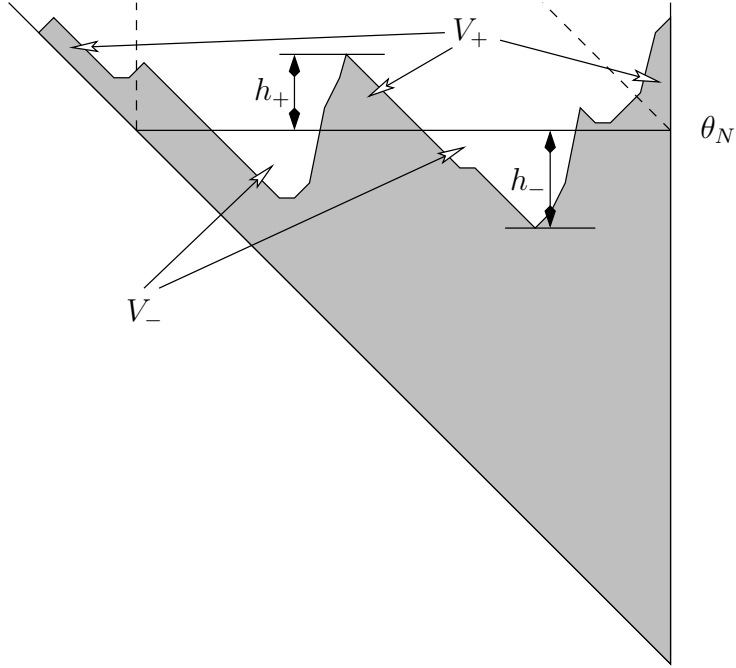
\begin{figure}[tb]
\begin{center}
\input{p3.pstex_t}
\caption{On the definition of the quantities $V_{\pm}(S)$, $h_{\pm}(S)$}
\label{p3}
\end{center}
\end{figure}

Define
\begin{eqnarray*}
h_+(S) &=& \max\{i: \mbox{ there exists } j\in[i-\theta_{|S|},\theta_{|S|}]\\
 &&~~~~~~~~~~~~~~~~~~~~~ \mbox{ such that } \RR_{i,j}(S)=1\} - \theta_{|S|},\\
h_-(S) &=& \theta_{|S|} - \min\{i: \mbox{ there exists~$j$ such that }
            \RR_{i,j}(S)=0 \},
\end{eqnarray*}
and
\begin{eqnarray*}
V_+(S) &=& |\{(i,j)\in\ZZ : i>\theta_{|S|}, \RR_{i,j}(S)=1\}|,\\
V_-(S) &=& |\{(i,j)\in\ZZ : i\leq\theta_{|S|}, \RR_{i,j}(S)=0\}|
\end{eqnarray*}
(cf.~Figure~\ref{p3}).
In words,
\begin{itemize}
\item $h_-(S)$ is the maximal vertical distance between~$\theta_N$ and
the holes below~$\theta_N$;
\item $h_+(S)$ is the maximal vertical distance between~$\theta_N$ and
the particles above~$\theta_N$ which also lie inside the area indicated
by the dashed lines on Figure~\ref{p3};
\item $V_-$ is the total area covered by the holes below~$\theta_N$;
\item $V_+$ is the total area covered by the particles above~$\theta_N$.
\end{itemize}
Similarly to the energy~$E(S)$, all those quantities could be used to
measure the deviation of~$S$ from the ``almost triangular''
configuration~$\T_0^N$.
Consider also the normalized quantities
$\th_{\pm}(S)=|S|^{-1/2}h_{\pm}(S)$, $\tV_{\pm}(S)=|S|^{-1}V_{\pm}(S)$,
and $\tE_{\pm}(S)=|S|^{-3/2}E_{\pm}(S)$, $\tE(S)=|S|^{-3/2}E(S)$.

\begin{lmm}
\label{alphas}
For all~$S\in G^{|S|}_{\gamma_1,\gamma_2}$ there exist constants
$\alpha_i, i=1,\ldots,8$ (depending on~$\gamma_1, \gamma_2$) such that
\begin{eqnarray}
\alpha_1 \th_-^3(S) ~ \leq & \tE_-(S) &\leq ~\alpha_2\th_-^2(S),\label{a12}\\
\alpha_3 \th_+^3(S) ~\leq &\tE_+(S) &\leq ~\alpha_4\th_+(S),\label{a34}\\
\alpha_5 \th_-^2(S) ~\leq &\tV_-(S)& \leq ~\alpha_6\th_-(S),\label{a56}\\
\alpha_7 \th_+^2(S) ~\leq &\tV_+(S) &\leq ~\alpha_8\th_+(S).\label{a78}
\end{eqnarray}
\end{lmm}

\noindent
{\it Proof.}
It is elementary to obtain
the inequalities~(\ref{a12}) and~(\ref{a56}) from~(\ref{pust_tri}).
Analogously, to obtain~(\ref{a34}) and~(\ref{a78}),
one can use~(\ref{poln_tri})
and the fact that $S\in G^{|S|}_{\gamma_1,\gamma_2}$
together with the following observation. If~$\RR_{i,j}(S)=1$ for
some~$i>\theta_{|S|}$, then either $i-\theta_{|S|}\leq h_+(S)$ or
$\min\{j,i-j\}\leq h_+(S)$.
\qed

\medskip

Consider a configuration~$S$ such that~$|S|=N$. By definition of~$\theta_N$,
there exists a constant~$\hat\beta>0$ such that
\begin{equation}
\label{hatbeta}
0\leq V_+(S) - V_-(S) \leq \hat\beta \sq.
\end{equation}
Also, we will always tacitly assume that $V_-(S) \geq \hat\beta \sq$,
i.e., we will not consider configurations that are ``too close''
to the triangle.
In this case there are constants~$\beta_1,\beta_2 >0$ such that
\begin{equation}
\label{beta12}
\beta_1 < \frac{\tV_+(S)}{\tV_-(S)} < \beta_2
\end{equation}
(note also that if~$N$ is a triangular number, then
$\tV_+(S)/\tV_-(S)=1$ for any $S\in\X_N$).
Using~(\ref{a12}), (\ref{a56}), and~(\ref{beta12}),
we obtain
\[
 \tE_-(S)\leq\alpha_2\th_-^2(S)\leq\frac{\alpha_2}{\alpha_5}\tV_-(S)
 \leq\frac{\alpha_2}{\alpha_5\beta_1}\tV_+(S)
  \leq \frac{\alpha_2\alpha_8}{\alpha_5\beta_1}\th_+(S),
\]
and, by~(\ref{a34}), $\tE_+(S)\leq \alpha_4\th_+(S)$. This shows that
there is a constant~$C_1$ such that $\th_+(S)\geq C_1 \tE(S)$
for all~$S\in G^{|S|}_{\gamma_1,\gamma_2}$.
Analogously, we obtain that for some~$C_2,C_3$ it holds that
$\th_-(S)\geq C_2\tE_+^2(S)$ and $\th_-(S)\geq C_3\tE_-^{1/2}(S)$.
By Lemma~\ref{prop_E} (i) the quantity~$\tE$ is bounded on
$G^{|S|}_{\gamma_1,\gamma_2}$, so there is~$C_4$ such that
$\tE_-^{1/2}(S)\geq C_4 \tE_-^2(S)$, which implies that
$\th_-(S)\geq C_5\tE^2(S)$ for some~$C_5$. Finally, we use
Lemma~\ref{prop_E} (i) once again to obtain that
there exists~$\beta=\beta(\gamma_1,\gamma_2)$ such that
\begin{equation}
\label{minh}
\min\{\th_+(S),\th_-(S)\} \geq \beta\tE^2(S)
\end{equation}
when~$S\in G^{|S|}_{\gamma_1,\gamma_2}$.

\subsection{Proof of Theorem~\ref{det_BS}}
\label{s_proof_det}
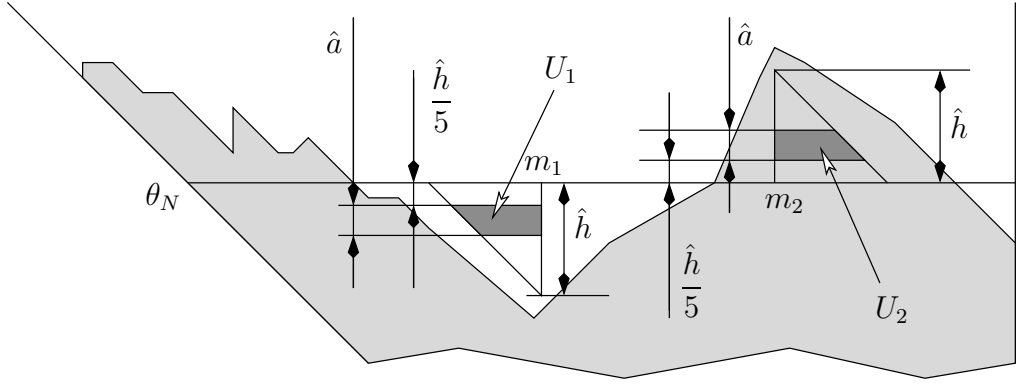
\begin{figure}
\begin{center}
\input{p4.pstex_t}
\caption{On the definition of sets $U_1$, $U_2$}
\label{p4}
\end{center}
\end{figure}

First, the idea is to prove
that after $O(\sq)$ moves, the ``normalized energy'' $\tE$ will decrease by a
considerable amount. Consider a configuration~$S$ with $|S|=N$.
Abbreviate $\hh = \lf \beta\tE^2(S)\sq\rf$; by~(\ref{minh}),
we can find~$m_1,m_2$ such that $\RR_{\theta_N-\hh,m_1}(S)=0$
and $\RR_{\theta_N+\hh,m_2}(S)=1$. Moreover, without loss
of generality one can suppose that~$\hh$ is divisible by~$5$.
Define also $j_1=m_1+4\hh/5$, $j_2=m_2-4\hh/5$,
and $\ha=\frac{\sqrt{2}\hh^2}{25\sq}$.
Define two sets $U_1,U_2\subset\ZZ$ by
\begin{eqnarray*}
U_1 &=& \Big\{(i,j) :
\theta_N-\ha-\frac{\hh}{5}\leq i \leq \theta_N-\frac{\hh}{5},
            m_1\leq j\leq j_1+i-\theta_N+\frac{\hh}{5} \Big\},\\
U_2 &=& \Big\{(i,j) :
\theta_N+\frac{\hh}{5}\leq i \leq \theta_N+\ha+\frac{\hh}{5},
             j_2+i-\theta_N-\frac{\hh}{5} \leq j \leq m_2\Big\},
\end{eqnarray*}
(see Figure~\ref{p4}).
Note that from~(\ref{pust_tri}) and~(\ref{poln_tri})
it follows that $\RR_{i,j}(S)=0$ for all $(i,j)\in U_1$ and
$\RR_{i,j}(S)=1$ for all $(i,j)\in U_2$.

Abbreviate also $i_0=\theta_N-\hh/5$, $i'_0=\theta_N+\hh/5$.
Now, the idea is to consider the evolution of
sets $U_1,U_2$ at times $i_0k$, $k=0,1,2,\ldots$. First, note that
${\hat J}_{i_0,j}(ki_0)=j$, for any~$j$ and~$k$. Then, each time we
make a complete turn (i.e., $i_0$ moves) a particle which was on the
level $\theta_N+\hh/5$ will be $2\hh/5$ units to the left of its initial
position (provided it did not fall). This shows that there exists
$k_0 \leq \frac{5\sqrt{2}}{2\hh}\sq$ such that
\[
  \big| [{\hat J}_{i'_0,j_2}(k_0i_0),{\hat J}_{i'_0,m_2}(k_0i_0)]
      \cap [m_1,j_1] \big|   \geq   \frac{2\hh}{5}
\]
(when ${\hat J}_{i'_0,j_2}(k_0i_0)>{\hat J}_{i'_0,m_2}(k_0i_0)$,
by $[{\hat J}_{i'_0,j_2}(k_0i_0),{\hat J}_{i'_0,m_2}(k_0i_0)]$
we mean in fact $[0,{\hat J}_{i'_0,m_2}(k_0i_0)]
      \cup[{\hat J}_{i'_0,j_2}(k_0i_0),i_0]$).
Take $j_3,j_4$ such that
\[
[j_3,j_4]\subset
\big([{\hat J}_{i'_0,j_2}(k_0i_0),{\hat J}_{i'_0,m_2}(k_0i_0)] \cap [m_1,j_1]\big)
\]
and $j_4-j_3+1 = \frac{2\hh}{5}$. We consider two cases:

\medskip
\noindent
\underline{Case 1}: at time~$k_0i_0$ in the set
\[
  U'_1 = \{(i,j): i\in [i_0-\ha,i_0], j\in[j_3,j_3-1+\hh/5]\}
\]
there is at least one hole, i.e.,
 $\RR_{i,j}(Q_1^{(k_0i_0)}S)=0$ for at least one $(i,j)\in U'_1$.
In this case, by~(\ref{pust_tri})
no particle can be in the set
\[
  U'_2 = \{(i,j): i\in [i'_0,i'_0+\ha], j\in[j_4-\hh/5,j_4]\},
\]
i.e., for all $(i,j)\in U'_2$ we have that $\RR_{i,j}(Q_1^{(k_0i_0)}S)=0$.
Note that $\ha k_0 \leq \hh/5$, so the ``image'' of~$U_2$ after~$k_0$
turns completely covers~$U'_2$. On the other hand, $U'_2$ must be
completely empty, so there should have been a lot of particle falls
in order to avoid~$U'_2$. In what follows we estimate
the minimal number of falls necessary (and, consequently, we find the
minimal amount by which the energy~$E$ should decrease).
Define the set
\[
  U''_2 = \{(i,j): i\in [i'_0+\ha/2,i'_0+\ha], j\in[j_4-\hh/10,j_4]\}
               \subset U'_2.
\]
For any $(i,j)\in U''_2$ there is a unique~$j'$ such that
${\hat J}_{i,j'}(k_0i_0)=j$, and, by the above observation,
$(i,j')\in U_2$, so the cell $(i,j')$ originally contained a particle.
To guarantee that that particle is not in~$U''_2$ at time~$k_0i_0$,
at least one of the following two possibilities must occur:
\begin{itemize}
\item either $\MM_{i,j'}(k_0i_0) \geq \ha/2$,
\item or ${\hat J}_{i,j'}(k_0i_0)-J_{i,j'}(k_0i_0)>\hh/10$, but in this case,
by Lemma~\ref{JM},
$\MM_{i,j'}(k_0i_0) \geq \frac{\hh}{10k_0}>
          \frac{\hh^2}{25\sqrt{2}\sq}=\ha/2$.
\end{itemize}
Denote $h_0=\hh N^{-1/2}$; for the both of the above possibilities,
we obtained in fact that $\MM_{i,j'}(k_0i_0) \geq C_1 h_0^2 \sq$.
Since the number of cells in the set~$U''_2$ is at least $C_2 h_0^3 N$,
the number of particle falls until time~$k_0i_0$ should be
at least $C_1 h_0^2 \sq \times C_2 h_0^3 N = C_1 C_2 h_0^5 N^{3/2}$.
By Lemma~\ref{prop_E} (ii), it means that, for the Case~1,
\begin{equation}
\label{hjk}
  \tE(Q_1^{(k_0i_0)}S) - \tE(S) \leq - C_1 C_2 h_0^5.
\end{equation}

\medskip
\noindent
\underline{Case 2}: there are no holes at time~$k_0i_0$ in the set~$U'_1$,
i.e., $\RR_{i,j}(Q_1^{(k_0i_0)}S)=1$ for all $(i,j)\in U'_1$.
Using the duality between holes and particles, this case can be treated
quite analogously to the Case~1. Namely, we note first that
the ``image'' of~$U_1$ after~$k_0$ turns completely covers~$U'_1$.
So, in order to escape~$U'_1$, the holes that are ``candidates'' to
be there must make a sufficient number of movements in the upwards direction.
In the same way
as in the Case~1, one can work out all the details to obtain that~(\ref{hjk})
is valid for the Case~2 as well.

\medskip

We continue proving Theorem~\ref{det_BS}.
By~(\ref{minh}) and~(\ref{hjk}), there exist~$\lambda_1,\lambda_2$
such that
\begin{equation}
\label{iteration}
  \tE(Q_1^{(n_S)}S) - \tE(S) \leq - \lambda_1 \tE^{10}(S),
\end{equation}
where $n_S=\lambda_2\tE^{-2}(S)\sq$ (the formula~(\ref{iteration}) will play an
important role in the proof of Theorem~\ref{ran_BS} as well). Consider now the initial
configuration~$S_0\in G^N_{\gamma_1,\gamma_2}$; by Lemma~\ref{prop_E} (i),
$a_0:=\tE(S_0)\leq \Gamma$ for some $\Gamma$. Fix an arbitrary~$\eps>0$ and
define $\phi(x) = x - \lambda_1 x^{10}$; then there exists~$k_1$
(depending only on~$\Gamma,\eps,\lambda_1$) such that $\phi^{(k_1)}(a_0)<\eps$.
By~(\ref{iteration}) this means that
\begin{equation}
\label{small_energy}
  \tE(Q_1^{(n'_S)}S) \leq \eps,
\end{equation}
where $n'_S = k_1\lambda_2\eps^{-2}\sq$, i.e., after~$O(\sq)$ moves
we will arrive to a confi\-gu\-ra\-tion with small normalized energy~$\tE$.

Now we are almost done with the proof of Theorem~\ref{det_BS},
and it remains only to make one small effort: we have to prove that if
the energy $\tE(S)$ is small, then either~$S$ is already close
to the triangular configuration~$\T(1,|S|)$ (in the sense of the
distance $\rho$), or it will come close to $\T(1,|S|)$ after
$O(\sq)$ moves.

Define the sets
\begin{eqnarray*}
\VV(\eps,N) &=& \big\{S : |S|=N,  \max\{h_+(S),h_-(S)\} \leq \eps\sq\big\},\\
{\hat \VV}(\eps,N) &=&\big\{S: |S|=N,
       \max\{i: \mbox{ there exists $j$ such that } \\
 &&      ~~~~~~~~~~~~~~~\RR_{i,j}(S)=1\}
            \leq \theta_N + \eps\sq \big\} \cap \VV(\eps,N).
\end{eqnarray*}
It is elementary to see that, for fixed~$\eps$ and for all~$N$ large enough
\begin{equation}
\label{subset}
{\hat \VV}(\eps,N) \subset \T(2\eps,1,N).
\end{equation}

We need the following
\begin{lmm}
\label{srezaem_hvosty}
Suppose that $S\in\VV(\eps,N)\cap G^N_{\gamma_1,\gamma_2}$, and put
$n_0 = [\max\{\gamma_1,\gamma_2\}]+2$.
Then $Q_1^{(n_0\theta_N)}S\in {\hat \VV}(2n_0\eps,N)$.
\end{lmm}

\noindent
{\it Proof.} Define
\begin{eqnarray*}
\VV'(\eps,N) &=&\big\{S: |S|=N,
       \max\{i: \mbox{ there exists $j\in[1,\theta_N]$ such that } \\
 &&      ~~~~~~~~~~~~~~~~~~~~~~~~~\RR_{i,j}(S)=1\}
            \leq \theta_N + \eps\sq \big\} \cap \VV(\eps,N).
\end{eqnarray*}
First, if $\eps<1/\sqrt{2}$ and $S\in\VV(\eps,N)\cap G^N_{\gamma_1,\gamma_2}$,
then the set $\{(i,j): i\geq\theta_N + \eps\sq,
 j\in [1,\eps\sq] \}$ will be empty of particles after~$\eps\sq$ moves. By
examining where those particles could go, we see that
$Q_1^{(\eps\sq)}S\in \VV'(2\eps,N)$ and that
\[
\max\{i: \mbox{ there exists $j$ such that }\RR_{i,j}(S')=1\}-\theta_N
       \leq (\eps+\max\{\gamma_1,\gamma_2\})\sq,
\]
where $S'=Q_1^{(\eps\sq)}S$.
To conclude the proof of the Lemma~\ref{srezaem_hvosty},
note the following two facts:
\begin{itemize}
\item Suppose that at some moment the
configuration belongs to the set $\VV'(\eps',N)$.
Then if $i'>\theta_N+\eps'\sq$
and there are some particles in the set $\{(i,j): i\geq i', j=i'\}$, then at
the next moment all those particles will fall one unit.
\item If $S\in\VV'(\eps',N)$ and
\[
\max_{i>\theta_N+\eps'\sq}(i- \min\{j: \RR_{i,j}(S)=1\}) \leq \eps''\sq,
\]
then $Q_1^{(\theta_N+\eps' \sq)}S\in \VV'(\eps'+\eps'',N)$
(to see this, it is enough to figure out what happens with the
configuration~$S''$ after $\theta_N+\eps' \sq$ moves,
where~$S''$ is defined by $\RR_{i,j}(S'')=1$ whenever either
$i\leq \theta_N+\eps'\sq$ or $i-j<\eps''\sq$).
\end{itemize}
\qed

\medskip

Now we are ready to finish the proof of Theorem~\ref{det_BS}.
  From~(\ref{a12}), (\ref{a34}), and~(\ref{small_energy}) we
obtain that, if the initial configuration belongs
to~$G^N_{\gamma_1,\gamma_2}$, then after~$O(\sq)$ moves it will
be in $\VV(\eps',N)\cap G^N_{\gamma_1+\sqrt{2},\gamma_2+\sqrt{2}}$,
where $\eps'=\eps^{1/3}/\min\{\alpha_1^{1/3},\alpha_3^{1/3}\}$
($\eps$ is from~(\ref{small_energy}), $\alpha_1,\alpha_3$ from Lemma~\ref{alphas}).
Applying Lemma~\ref{srezaem_hvosty},
we conclude the proof of Theorem~\ref{det_BS}.
\qed

\subsection{Proof of Theorem~\ref{ran_BS}}
\label{s_proof_rand}
Consider a finite irreducible discrete-time Markov chain with state
space~$X$, transition matrix~$P$, and stationary measure~$\pi$.
The following elementary fact will be useful in the course of the proof
of Theorem~\ref{ran_BS}: for any $A\subset X$ and $n\geq 1$
\begin{equation}
\label{osn_fakt}
\sum_{\substack{x\in A,\\ y\in A^c}}\pi(x)P_{xy}^{(n)} =
 \sum_{\substack{x\in A,\\ y\in A^c}}\pi(y)P_{yx}^{(n)}.
\end{equation}

Let us describe the main steps of the proof of Theorem~\ref{ran_BS}:
\begin{itemize}
\item first, in Lemma~\ref{reasonable} we prove (using (\ref{osn_fakt}))
that a typical configuration of the random game should be reasonable,
i.e., the number of piles and the biggest pile should be $O(\sq)$;
\item then, the idea is the following: starting from a reasonable configuration,
the macroscopic evolution of the profiles will be very similar in the
random game and in the deterministic game where the initial sizes of
the piles are~$p^{-1}$ times bigger (indeed, if the initial size of the pile in
the random game is $k=O(\sq)$, then it will be emptied typically in time
$k/p \pm O(N^{1/4})$);
\item unfortunately, it seems to be difficult to dominate the stochastic game
by the deterministic one directly. So, we introduce another deterministic
process by allowing the immigration of particles to the system on each step.
In Lemma~\ref{dominacao} we prove that the random Bulgarian solitaire is
in some sense dominated by this new deterministic process;
\item it is then possible to see that the process with immigration of
particles does not differ much (when the time interval in question
is not too long) from the deterministic Bulgarian solitaire,
because the total number of added particles is relatively small,
and they cannot be very concentrated (Lemma~\ref{nao_fugira} takes care of
the latter statement). Using this observation and applying
inequality~(\ref{iteration}) from the previous section, we obtain that
the (suitably defined) energy will decrease with large probability
after a sufficiently large number of steps (this is Lemma~\ref{iter_rand});
\item the rest of the proof is a straightforward (although somewhat
lengthy) application of~(\ref{osn_fakt}) and Lemma~\ref{iter_rand}.
\end{itemize}

So, the first step is to prove that a typical configuration $S\in\X_N$
should be ``reasonable'', i.e., it should belong to $G^N_{\gamma'_1,\gamma'_2}$
for some $\gamma'_1,\gamma'_2$:
\begin{lmm}
\label{reasonable}
\begin{itemize}
\item[(i)]
For any~$p\in (0,1)$ there exist positive constants $\sigma_0,\gamma'_1,\gamma'_2$
(depending on~$p$) such that
\begin{equation}
\label{sdf}
 \pi_{p,N} (G^N_{\gamma'_1,\gamma'_2}) \geq 1- e^{-\sigma_0\sq}.
\end{equation}
for all~$N$.
\item[(ii)] Also, suppose that $S\in G^N_{\gamma'_1,\gamma'_2}$, where
$\gamma'_1,\gamma'_2$ are the quantities from item~(i). Then
there exist $\gamma''_1,\gamma''_2$ and $\sigma_1$ such that
for any $M>1$
\[
 \Pr[Q_p^{(n)}S \in G^N_{\gamma''_1,\gamma''_2}] \geq 1-N^M e^{-\sigma_1\sq}.
\]
\end{itemize}
\end{lmm}

\noindent
{\it Proof.}
We begin by proving~(i). Consider the sets
\begin{eqnarray*}
 A_N &=& \{S: |S|=N, \ell(S)>3\sq/p\}\\
 A'_N &=& \{S: |S|=N, \ell(S)>5\sq/p\}.
\end{eqnarray*}
Note that if $S=(k_1,\ldots,k_{\ell(S)})\in A_N$, then
$|\{i:k_i>p\sq\}|<\sq/p$, i.e., in~$S$ there are at most~$\sq/p$
piles with at least~$p\sq$ cards. Clearly, if a pile had no more than~$p\sq$
cards, then there is~$C_1>0$ such that by the time~$3\sq/2$
that pile will be empty with
probability at least $1-e^{-C_1\sq}$.
During the time~$3\sq/2$ only~$3\sq/2$ new piles can appear,
so, since $1/p+3/2<3/p$, for any~$S\in A_N$ we have that,
\begin{equation}
\label{A1}
\Pr[Q_p^{(3\sq/2)}S\in \X_N\setminus A_N] \geq 1 - N e^{-C_1\sq}.
\end{equation}
Since $A'_N\subset A_N$, (\ref{A1}) also shows that for
any~$S\in A_N\setminus A'_N$
we have $\Pr[Q_p^{(3\sq/2)}S\in A'_N] \leq N e^{-C_1\sq}$. On the other hand,
if $S\in \X_N\setminus A_N$, then clearly $\Pr[Q_p^{(3\sq/2)}S\in A'_N]=0$,
so (since $3/2<2/p$) for any $S\in \X_N\setminus A'_N$ we have
\begin{equation}
\label{A2}
\Pr[Q_p^{(3\sq/2)}S\in A'_N] \leq N e^{-C_1\sq}.
\end{equation}
Now we use~(\ref{osn_fakt}) with $A=A'_N$ and $n=3\sq/2$
to obtain from~(\ref{A1})
and~(\ref{A2}) that for some $C_2>0$
\begin{equation}
\label{AAA}
 \pi_{p,N}(\X_N\setminus A'_N) = \pi_{p,N}(S: \ell(S)\leq 5\sq/p)
          \geq 1 - e^{-C_2\sq}.
\end{equation}

For $k=1,\ldots, \lfloor\sq\rfloor$ define
\[
 B_N^{(k)} = \{S: |S|=N, (k-1)p\sq < R(S) \leq kp\sq\}\cap (\X_N\setminus A'_N).
\]
Suppose that $S\in B_N^{(k)}$ for some~$k>1+\frac{5}{p^2}+\frac{3}{2p}$,
and let us try to figure out what the configuration $Q_p^{(3\sq/2)}S$
should look like. Note that
\begin{itemize}
\item $\ell(S)\leq \frac{5}{p}\sq$, and moreover
$\ell(Q_p^{(m)}S)\leq (\frac{5}{p}+\frac{3}{2})\sq$
for all $m\leq\frac{3}{2}\sq$,
so, since $(k-1)p>\frac{5}{p}+\frac{3}{2}$, no new pile of size greater than
$(k-1)p\sq$ can appear until the moment $\frac{3}{2}\sq$;
\item the evolution of a single pile can be modeled by a random walk
on~$\Z_+$ which jumps one unit to the left with probability~$p$ and
holds its position with probability~$1-p$. This shows that if the size
of the pile was less than $kp\sq$, then after $\frac{3}{2}\sq$ moves
it will be less than $(k-1)p\sq$ with probability at least
$1 - e^{-C_3\sq}$ for some $C_3>0$.
\end{itemize}
  From the above facts we deduce that for any $S\in B_N^{(k)}$
\[
  \Pr[Q_p^{(3\sq/2)}S\notin B_N^{(k)}] \leq N e^{-C_3\sq},
\]
as long as $k>1+\frac{5}{p^2}+\frac{3}{2p}$.
Now using~(\ref{osn_fakt}) with $A=B_N^{(k)}$ and $n=3\sq/2$, we
obtain that for some $C_4>0$
\[
 \pi_{p,N}(B_N^{(k)}) \leq e^{-C_4\sq} + \sum_{m\geq k}\pi_{p,N}(B_N^{(m)}),
\]
so by induction one can show that
$\pi_{p,N}(B_N^{(k)}) \leq C_5 N^2 e^{-C_4\sq}$.
Summing over $k>1+\frac{5}{p^2}+\frac{3}{2p}$ and recalling~(\ref{AAA}),
we conclude the proof of the part~(i) of Lemma~\ref{reasonable}
(with $\gamma'_1=\frac{5}{p}$, $\gamma'_2=\frac{5}{p}+\frac{3}{2}$).

To prove the part~(ii), first observe that in the proof of~(i)
we have constructed $\gamma'_1,\gamma'_2$ in such a way that
for any $S\in G^N_{\gamma'_1,\gamma'_2}$
\[
  \Pr[Q_p^{(3\sq/2)}S \in \X_N\setminus G^N_{\gamma'_1,\gamma'_2}]
              \leq e^{-C\sq}.
\]
To complete the proof of~(ii), it is enough to take $\gamma''_1=\gamma'_1+\frac{3}{2}$,
$\gamma''_2=\max\{\gamma'_1,\gamma'_2\}+\frac{3}{2}$ (note that for any~$S$
we have $\ell(Q_pS)-\ell(S)\leq 1$, $R(Q_pS)\leq\max\{R(S),\ell(S)\}$).
\qed

\medskip

We continue proving Theorem~\ref{ran_BS}. Now, the main idea is
the following: first,
to dominate the random Bulgarian solitaire by a certain deterministic
process (that we will call Bulgarian solitaire with immigration
of particles), and then apply to that process some methods from the
proof of Theorem~\ref{det_BS}.

Fix $\delta_0 \in (0,\frac{1}{4}-36a)$, and abbreviate
$\kk=\lc N^{\delta_0+\frac{1}{4}}\rc$.
Denote also~$N_p:=[N/p]$.
For any~$S$, let us define
configurations~$\D(S),\tD(S)$ in the following way:
if $S=(n_1,\ldots,n_{\ell(S)})$, let
\begin{eqnarray*}
\D(S) &=& \Big(\Big\lc\frac{n_1}{p}\Big\rc-z_1,\ldots,
                      \Big\lc\frac{n_{\ell(S)}}{p}\Big\rc-z_{\ell(S)}\Big),\\
\tD(S) &=& \Big(\Big\lc\frac{n_1}{p}\Big\rc-z_{1}+\kk,\ldots,
                      \Big\lc\frac{n_{\ell(S)}}{p}\Big\rc-z_{\ell(S)}+\kk\Big),
\end{eqnarray*}
where $z_i=z_i(S)\in\{0,1\}$ are chosen in such a way that
$z_1\geq z_2 \geq \ldots \geq z_{\ell(S)}$ and for any
$S\in\X_N$ we have $|\D(S)|=N_p$.
Define the operator~$\tQ$ by
\[
 \tQ S = \ord(n_1-1,\ldots,n_{\ell(S)}-1,\ell(S)+\kk),
\]
i.e., making the $\tQ$-move consists of making a move of deterministic
Bulgarian solitaire, and then adding~$\kk$ cards to the new pile
(so that $|\tQ S|-|S|=\kk$). For the simplicity of notations, we do not
indicate in~$\tQ$ the dependence on~$N$ and $\delta_0$; note also that in
the above display we do not assume that $|S|=N$, so~$\tQ$ need not
apply to only $S\in\X_N$.

 For two configurations~$S_1=(n_1,\ldots,n_{\ell(S_1)})$,
$S_2=(m_1,\ldots,m_{\ell(S_2)})$ we say that $S_1\leq S_2$ if
$\ell(S_1)\leq \ell(S_2)$ and $n_j\leq m_j$ for all $j=1,\ldots,\ell(S_1)$.

\begin{lmm}
\label{dominacao}
Suppose that $|S|=N$ and $S\in G^N_{\gamma'_1,\gamma'_2}$
 (where $\gamma'_1,\gamma'_2$
are the quantities from Lemma~\ref{reasonable}).
Then for any $M>0$ there exists~$\sigma_2$
such that
\begin{equation}
\label{nkr}
\Pr[\tQ^{(n)}\tD(S)\geq \D(Q_p^{(n)}S) \mbox{ for all }n\leq N^M] \geq
      1 - N^M e^{-\sigma_2N^{\delta_0}}.
\end{equation}
\end{lmm}

\noindent
{\it Proof.}
Let us refer to the $\ell(S)$ piles of~$S$ and $\tD(S)$ as
$\PP_1,\ldots,\PP_{\ell(S)}$ and $\tPP_1,\ldots,\tPP_{\ell(S)}$
respectively. Then, the piles born at the moment~$i$ are referred to
as $\PP_{\ell(S)+i}$ and $\tPP_{\ell(S)+i}$.
Using the notation $(x)^+:=\max\{x,0\}$, for $n\geq (i-\ell(S))^+$,
let $\PP_i(n)$ and $\tPP_i(n)$ stand for the sizes of the piles
 $\PP_i$ and $\tPP_i$ at the moment~$n$, respectively (if a pile is
emptied at some moment $n^*<n$, then we mean that the size remains~$0$
for all $m\geq n^*$).

Clearly, the event
\begin{equation}
\label{mnbvs}
\{\tQ^{(n)}\tD(S)\geq \D(Q_p^{(n)}S) \mbox{ for all }n\leq N^M\}
      \subset \bigcap_{i=1}^{N^M+\ell(S)} \Lambda_i,
\end{equation}
where we define the event~$\Lambda_i$ by
\[
 \Lambda_i = \{\PP_i(n) \leq p\tPP_i(n)\mbox{ for all } n\geq (i-\ell(S))^+\}.
\]
Define also the event $D=\{Q_p^{(n)}S\in G^N_{\gamma''_1,\gamma''_2}\}$;
by Lemma~\ref{reasonable} (ii) we know that
$\Pr[D]\geq 1 - N^M e^{-\sigma_1\sq}$. On the other hand,
\begin{equation}
\label{uyi**}
\Pr[\Lambda_i\mid \Lambda_1,\ldots,\Lambda_{i-1},D] \geq
\Pr[H_i]\Pr[\Lambda_i\mid H_i,D],
\end{equation}
where $H_i=\{p\tPP_i((i-\ell(S))^+) \geq \PP_i((i-\ell(S))^+)+\kk/2\}$.
Now, on~$D$ we have that $\ell(Q_p^{(n)}S)=O(\sq)$ for all $n\leq N^M$,
and on $\Lambda_1\cap\ldots\cap\Lambda_{i-1}$ it holds that
$\ell(\tQ^{(n)}\tD(S))\geq\ell(Q_p^{(n)}S)$.
Using the Large Deviation bound for the Binomial distribution, we see
that the first term in the right-hand side of~(\ref{uyi**})
is at least $1-e^{-C_1N^{\delta_0}}$. As for the second term, note that
the difference between $p\tPP_i(\cdot)$ and $\PP_i(\cdot)$ is
a random walk with drift~$0$. Since the time that the pile $\PP_i$
needs to be emptied is $O(\sq)$, the second term in~(\ref{uyi**})
is in fact the probability that such a random walk does not deviate
 from its initial position by more than~$\kk/2$ by time $O(\sq)$;
clearly, that probability is bounded from below by
$1-e^{-C_2N^{\delta_0}}$. Then, it is immediate to deduce
Lemma~\ref{dominacao} from~(\ref{mnbvs}) and~(\ref{uyi**}).
\qed

\medskip
Recall that (cf.\ the proof of Lemma~\ref{prop_E})
for any configuration~$S$ we use the notation
\[
 \H(S) = \max\{i: \mbox{ there exists $j$ such that } \RR_{i,j}(S)=1\}.
\]
\begin{lmm}
\label{nao_fugira}
Suppose that $S\in G^N_{\gamma'_1,\gamma'_2}$ and let~${\tilde\beta}$ be
such that ${\tilde\beta} \leq \frac{1}{4}-\delta_0$. Then there exists~$L_0$
 such that $\H(\tQ^{(n)}\tD(S))\leq L_0\sq$
for all $n\leq N^{\frac{1}{2}+{\tilde\beta}}$.
\end{lmm}

\noindent
{\it Proof.} Let $b_0=\H(\tD(S))$ and denote $\hat b_0 = b_0(b_0+1)/2$.
Define the triangular configuration~$\hT_0$ by
$\RR_{i,j}(\hT_0)=\1\{i\leq \hat b_0\}$; then, clearly,
$\tD(S)\leq \hT_0$. Denote $b_1 = b_0 + \lc\kk/\sqrt{2}\rc$,
$\hat b_1 = b_1(b_1+1)/2$ and define the configuration~$\hT_1$ by
$\RR_{i,j}(\hT_1)=\1\{i\leq \hat b_1\}$. By examining the
$\tQ$-evolution of~$\hT_0$ on the Etienne diagram, it is clear
that $\tQ^{(n)}\hT_0\leq \hT_1$ for all $n\leq b_0$. We then
repeat this construction by defining $b_{m+1} = b_m + \lc\kk/\sqrt{2}\rc$,
$\hat b_{m+1} = b_{m+1}(b_{m+1}+1)/2$ and the configuration~$\hT_{m+1}$ by
$\RR_{i,j}(\hT_{m+1})=\1\{i\leq \hat b_{m+1}\}$. Analogously, we
obtain that $\tQ^{(n)}\hT_m\leq \hT_{m+1}$ for all $n\leq b_n$.
A simple monotonicity argument then shows that
$\tQ^{(n)}\tD(S)\leq \hT_{m+1}$ for all $n\leq b_0+\cdots+b_m$.
We have $b_0+\cdots+b_m\geq (m+1)b_0$ and $b_0\leq C_1\sq$ for
some~$C_1$, so $\tQ^{(n)}\tD(S)\leq \hT_{C_1^{-1}N^{\tilde\beta}}$ for all
$n\leq N^{1/2+{\tilde\beta}}$. So,
since $\frac{1}{4}+\delta_0+{\tilde\beta} \leq \frac{1}{2}$, for some~$L_0$
we have
\[
 \H(\tQ^{(n)}\tD(S)) \leq C_1\sq +
         \Big\lc\frac{\kk}{\sqrt{2}}\Big\rc C_1^{-1}N^{\tilde\beta}
              \leq L_0\sq
\]
for all~$N$, thus concluding the proof of Lemma~\ref{nao_fugira}.
\qed

\begin{lmm}
\label{iter_rand}
Fix some~$\tilde a\in (0,1/16)$ and suppose that a
configuration~$S\in G^N_{\gamma'_1,\gamma'_2}$ is
such that $\tE(\D(S)) \geq \lambda_3N^{-\tilde a}$ for some~$\lambda_3$.
Then, with $\lambda_1,\lambda_2$ from~(\ref{iteration}),
we have for some~$\sigma_3,\delta_1>0$
\begin{equation}
\label{eq_l2.8}
\Pr[\tE(\D(Q_p^{(n'_S)}S))-\tE(\D(S)) \leq -\lambda_1\tE^{10}(\D(S))/2]
         \geq 1 - e^{-\sigma_3N^{\delta_1}},
\end{equation}
where $n'_S=\lambda_2\tE^{-2}(\D(S))\sq/p$. Moreover, (\ref{eq_l2.8})
remains true when $Q_p^{(n'_S)}$ is substituted by $Q_p^{(n)}$,
for any $n\in [n'_S,2n'_S]$.
\end{lmm}

\noindent
{\it Proof.}
First, each particle added to~$\D(S)$ changes the energy~$E$ by
at most~$O(\sq)$, so we have for some constants $C_1,C_2$
\begin{equation}
\label{asd}
|E(\tD(S))-E(\D(S))| \leq \ell(\D(S))\times C_1\sq\kk
  \leq C_2 N^{5/4+\delta_0}.
\end{equation}
Using the same sort of argument and the fact that
$\tQ^{(m)}S' \geq Q_1^{(m)}S''$ for any $m, S'\geq S''$, with the
help of Lemma~\ref{nao_fugira} and~(\ref{asd}) we obtain
\begin{equation}
\label{asd'}
|E(Q_1^{(n'_S)}\D(S))-E(\tQ^{(n'_S)}\tD(S))|
       \leq C_2 N^{5/4+\delta_0} + C_3 n'_S N^{3/4+\delta_0}.
\end{equation}
Introduce the event $D_1 = \{\tQ^{(n'_S)}\D(S) \geq \D(Q_p^{(n'_S)}S)\}$.
By Lemma~\ref{dominacao} we have
\begin{equation}
\label{fff***}
\Pr[D_1] \geq 1-n'_S e^{-\sigma_2N^{\delta_0}},
\end{equation}
and, since
\[
 |\tQ^{(n'_S)}\tD(S)|-|\D(Q_p^{(n'_S)}S)| \leq (O(\sq) + n'_S)\kk,
\]
analogously to~(\ref{asd})--(\ref{asd'}) we obtain that on~$D_1$
\begin{equation}
\label{fff****}
 |E(\tQ^{(n'_S)}\tD(S)) - E(\D(Q_p^{(n'_S)}S))|
           \leq C_4 n'_S N^{\frac{3}{4}+\delta_0}.
\end{equation}

Now, we have that $\tE^{10}(\D(S)) \geq \lambda_3^{10}N^{-10{\tilde a}}$, and
$n'_SN^{\frac{3}{4}+\delta_0}\leq C_5 N^{\frac{5}{4}+2{\tilde a}+\delta_0}$.
Since $\frac{3}{2}-10{\tilde a}>\frac{5}{4}+2{\tilde a}+\delta_0$, we
obtain the proof of~(\ref{eq_l2.8}) from~(\ref{iteration}), (\ref{asd'}),
(\ref{fff***}), and~(\ref{fff****}).

As for the second claim of Lemma~\ref{iter_rand},
we note that for $n\geq n'_S$,
by Lemma~\ref{prop_E} (ii) it holds that
 $\tE(Q_1^{(n)}\D(S))\leq \tE(Q_1^{(n'_S)}\D(S))$,
and then use the same kind of estimates
as used above.
\qed

\medskip

Now we are ready to finish the proof of Theorem~\ref{ran_BS}.
By Lemma~\ref{reasonable} (i) there are $\sigma_0,\gamma'_1,\gamma'_2$
such that~(\ref{sdf}) holds. Note that there exists
$\Gamma'=\Gamma'(\gamma'_1,\gamma'_2)$ such that if
$S\in G^N_{\gamma'_1,\gamma'_2}$,
then $\tE(\D(S))\leq \Gamma'$. Define $\psi(x)=x-\frac{1}{2}\lambda_1x^{-10}$.
Let $y_0=\Gamma'$ and $y_{i+1}=\psi(y_i)$ for $i\geq 0$.
Take $\eps=N^{-a}$, $a<1/144$, and define
$\hat\eps=\min\{\alpha_1,\alpha_3\}\eps^3$ (cf.~(\ref{a12}) and~(\ref{a34})).
Let $\hat n=\min\{n:y_n<\hat\eps\}$;
since $\hat\eps=O(N^{1/48})$, by examining the iteration scheme
$x\mapsto \psi(x)$ we obtain that
there exists $C_1$ such that $\hat n \leq C_1 N^{5/24}$.
Let
\[
\L_n = \{S\in G^N_{\gamma'_1,\gamma'_2} : \tE(\D(S))\in(y_{n+1},y_n]\},
\]
and define also $\L_{>n} = \displaystyle\bigcup_{k>n} \L_k$,
$\L_{<n} = \displaystyle\bigcup_{k<n} \L_k$.
Take any $n\leq\hat n$ and denote $m_n=\lambda_2y_{n+1}^{-2}\sq/p$.
By~(\ref{osn_fakt}) and Lemma~\ref{reasonable} (i) we can write
\begin{eqnarray}
\sum_{\substack{S_1\in \L_n\\ S_2\in \L_n^c}}
\pi_{p,N}(S_1)P_{S_1S_2}^{(m_n)} &=&
 \sum_{\substack{S_1\in \L_n\\ S_2\in \L_n^c}} \pi_{p,N}(S_2)P_{S_2S_1}^{(m_n)}
 \label{linha1}\\
 &\leq& e^{-\sigma_0\sq} + T_1 + T_2,
 \label{linha2}
\end{eqnarray}
where
\begin{eqnarray*}
T_1 &=& \sum_{\substack{S_1\in \L_n\\ S_2\in \L_{>n}}}
\pi_{p,N}(S_2)P_{S_2S_1}^{(m_n)},\\
T_2 &=& \sum_{\substack{S_1\in \L_n\\ S_2\in \L_{<n}}}
\pi_{p,N}(S_2)P_{S_2S_1}^{(m_n)}.
\end{eqnarray*}
Now, by Lemma~\ref{iter_rand}, the left-hand side of~(\ref{linha1})
can be bounded from below as follows:
\begin{eqnarray}
\sum_{\substack{S_1\in \L_n\\ S_2\in \L_n^c}}
\pi_{p,N}(S_1)P_{S_1S_2}^{(m_n)} &\geq&
 \sum_{S_1\in \L_n}\pi_{p,N}(S_1)\sum_{S_2\in \L_n^c}P_{S_1S_2}^{(m_n)}
\nonumber\\
 &=& \sum_{S_1\in \L_n}\pi_{p,N}(S_1)
        \Pr[Q_p^{(m_n)}S_1 \in \L_n^c]
\nonumber\\
 &\geq&\pi_{p,N}(\L_n) (1-e^{-\sigma_3N^{\delta_1}}). \label{yyy}
\end{eqnarray}
Again using Lemma~\ref{iter_rand}, we write
\begin{equation}
\label{yyy'}
 T_1 \leq \sum_{S_2\in \L_{>n}}\pi_{p,N}(S_2)e^{-\sigma_3N^{\delta_1}}
                  \leq e^{-\sigma_3N^{\delta_1}}.
\end{equation}
Using now~(\ref{yyy}) and~(\ref{yyy'}) together with the trivial bound
$T_2\leq \pi_{p,N}(\L_{<n})$, we obtain from~(\ref{linha1})--(\ref{linha2})
that for some $C_2>0$
\begin{equation}
\label{qqq}
 \pi_{p,N}(\L_n) \leq C_2 \Big(e^{-\sigma_0\sq} + e^{-\sigma_3N^{\delta_1}}
          + \sum_{k<n}\pi_{p,N}(\L_k)\Big).
\end{equation}
By induction, we then obtain that there is~$C_3>0$ such that for any $n<\hat n$
\begin{equation}
\label{qqq'}
\pi_{p,N}(\L_n) \leq C_3 n^2 e^{-\sigma_3N^{\delta_1}},
\end{equation}
so, since $\hat n \leq C_1 N^{5/24}$, taking summation
in~(\ref{qqq'}) we obtain for some~$C_4>0$ that
\begin{equation}
\label{oc_en}
\pi_{p,N}(S\in\X_N: \tE(\D(S))\geq\hat\eps) \leq C_4 N^{15/24} e^{-\sigma_3N^{\delta_1}}.
\end{equation}

Now, the last step of the proof of Theorem~\ref{ran_BS}
is analogous to what was done in Lemma~\ref{srezaem_hvosty}.
Note that if $\tE(\D(S))<\hat\eps$, then
\[
\max\{\th_-(\D(S)),\th_+(\D(S))\} \leq
         \Big(\frac{\hat\eps}{\min\{\alpha_1,\alpha_3\}}\Big)^{1/3} = \eps,
\]
so if $\tE(\D(S))<\hat\eps$, then $\D(S)\in\VV(\eps,|\D(S)|)$, thus showing that
\begin{equation}
\label{jsueh}
\pi_{p,N}(S\in\X_N: \D(S)\in\VV(\eps,N_p))
                    \geq 1 - C_4 N^{15/24} e^{-\sigma_3N^{\delta_1}}.
\end{equation}

Define
\begin{eqnarray*}
H_0 &=& \big\{S: \D(S)\in \VV(\eps,|\D(S)|), \max\{i:\mbox{ there
exists } j\leq \eps|\D(S)|^{1/2} \\
&& ~~~~~~~~~~~~~~~~ \mbox{ such that } \RR_{i,j}(\D(S))=1\} \geq
 \theta_{|\D(S)|}+2\eps|\D(S)|^{1/2}\big\}.
\end{eqnarray*}
Take any $S\in G^N_{\gamma'_1,\gamma'_2}$,
and denote $W=\{S \in \X_N:\D(S)\in\VV(\eps,N_p)\}$ (recall that $|\D(S)|=N_p$).
Using~(\ref{osn_fakt}), we write
\begin{eqnarray}
\sum_{\substack{S_1\in H_0\\ S_2\in H_0^c}}\pi_{p,N}(S_1)
       P_{S_1S_2}^{(2\eps N_p^{1/2})} &=&
\sum_{\substack{S_1\in H_0\\ S_2\in H_0^c}}\pi_{p,N}(S_2)
       P_{S_2S_1}^{(2\eps N_p^{1/2})} \label{lll1}\\
&=& T'_1 + T'_2, \label{lll2}
\end{eqnarray}
where
\begin{eqnarray*}
T'_1 &=& \sum_{\substack{S_1\in H_0\\ S_2\in H_0^c\cap W}}\pi_{p,N}(S_2)
       P_{S_2S_1}^{(2\eps N_p^{1/2})},\\
T'_2 &=& \sum_{\substack{S_1\in H_0\\ S_2\in H_0^c\cap W^c}}\pi_{p,N}(S_2)
       P_{S_2S_1}^{(2\eps N_p^{1/2})}.
\end{eqnarray*}
Observe that if $\D(S)\in \VV(\eps,N_p)$
and~$\eps$ is small enough, then after $2\eps N_p^{1/2}$ moves
there will be no particles in the set $\{(i,j): i\geq
\theta_{N_p}+2\eps N_p^{1/2}\}, j\leq \eps N_p^{1/2}$,
with probability at least $1-e^{-C_5\sq}$ for some~$C_5$.
So, for the left-hand side of~(\ref{lll1}) we can write
\begin{equation}
\label{lhs1}
\sum_{\substack{S_1\in H_0\\ S_2\in H_0^c}}\pi_{p,N}(S_1)
       P_{S_1S_2}^{(2\eps N_p^{1/2})} \geq \pi_{p,N}(H_0) (1-e^{-C_5\sq}).
\end{equation}
On the other hand, the same argument implies that $T'_1\leq e^{-C_5\sq}$
and the bound $T'_2 \leq \pi_{p,N_p}(W^c)$ is trivial. So, using~(\ref{oc_en})
and~(\ref{lhs1}), we obtain from~(\ref{lll1}) that
\begin{equation}
\label{oc_H0}
 \pi_{p,N}(H_0) \leq C_6 N^{15/24} e^{-\sigma_3N^{\delta_1}}.
\end{equation}

Abbreviate $\hat H = H_0^c\cap W$ and define
\begin{eqnarray*}
 F_k &=& \big\{S: \D(S)\in \hat H, \max\{i: \mbox{ there
   exists $j\geq \theta_{N_p}$ such that }\\
   && ~~~~~~~~~\RR_{i,j}(\D(S))=1\}-\theta_{N_p}-\eps\sq
        \in(2\eps k\sq, 2\eps(k+1)\sq]\big\},
\end{eqnarray*}
$F_{<k}=\displaystyle\bigcup_{m<k}F_k$, $F_{>k}=\displaystyle\bigcup_{m>k}F_k$.
Analogously to (\ref{linha1})--(\ref{linha2}) and (\ref{lll1})--(\ref{lll2}),
we write
\begin{eqnarray}
\sum_{\substack{S_1\in F_k\\ S_2\in F_k^c}}\pi_{p,N}(S_1)
       P_{S_1S_2}^{(4\eps N_p^{1/2})} &=&
\sum_{\substack{S_1\in F_k\\ S_2\in F_k^c}}\pi_{p,N}(S_2)
       P_{S_2S_1}^{(4\eps N_p^{1/2})} \label{zzz1}\\
&=& T''_1 + T''_2 + T''_3, \label{zzz2}
\end{eqnarray}
where
\begin{eqnarray*}
T''_1 &=& \sum_{\substack{S_1\in F_k\\ S_2\in F_{>k}}}\pi_{p,N}(S_2)
       P_{S_2S_1}^{(4\eps N_p^{1/2})},\\
T''_2 &=& \sum_{\substack{S_1\in F_k\\ S_2\in F_{<k}}}\pi_{p,N}(S_2)
       P_{S_2S_1}^{(4\eps N_p^{1/2})},\\
T''_3 &=& \sum_{\substack{S_1\in F_k\\ S_2\in {\hat H}^c}}\pi_{p,N}(S_2)
       P_{S_2S_1}^{(4\eps N_p^{1/2})}.
\end{eqnarray*}
The following fact can be deduced from~(\ref{poln_tri}): if $\D(S)\in\hat H$
and $\RR_{i,j}(\D(S))=1$ for some $i>\theta_{N_p}+\frac{5}{2}\eps N_p^{1/2}$,
then $i-j\leq \eps N_p^{1/2}$. Then,
by examining the evolution of $\D(S)$ on the Etienne diagram and
using Lemma~\ref{dominacao}, it is elementary to obtain that
for any $S\in F_k$, $k\geq 1$
\[
\Pr[\D(Q_p^{(4\eps N_p^{1/2})}S)\notin F_k \cup F_{>k}] \geq 1 - e^{-C_7\sq}.
\]
Using that fact, one can bound the left-hand side of~(\ref{zzz1}) from
below by $\pi_{p,N}(F_k) (1 - e^{-C_7\sq})$ and the term $T''_2$ can be
bounded from above by $e^{-C_7\sq}$. Then, it is straightforward to write
$T''_1\leq \sum_{m>k} \pi_{p,N}(F_m)$, $T''_3\leq\pi_{p,N}({\hat H}^c)$.
Denoting now $\tilde m = \frac{\gamma'_1\sq}{2\eps\sq}=\frac{\gamma'_1}{2\eps}$,
analogously to (\ref{qqq})--(\ref{qqq'}) we obtain
\begin{equation}
\label{ghhrt}
 \pi_{p,N}(F_k) \leq C_8 (\tilde m - k)^2 N^{15/24} e^{-\sigma_3N^{\delta_1}}.
\end{equation}
Summing over $k=1,\ldots,\tilde m$ and taking~(\ref{jsueh}) and~(\ref{oc_H0})
into account, we finally obtain that for some $C_9,\delta>0$ (depending on~$a$)
\[
 \pi_{p,N} (\hat\VV(3\eps,N_p)) \geq 1 - e^{-C_9N^\delta}.
\]
Since $\eps=N^{-a}$ and $a<1/144$ is arbitrary, we complete the proof
of Theorem~\ref{ran_BS}
(note that $\D^{-1}(\hat\VV(\eps,N_p))\subset\T(2\eps,p,N)$
for $\eps\gg N^{-1/2}$).
\qed

\section{Final remarks}
\label{s_fr}
A natural question that one may ask is: starting from
an initial configuration~$S$ with $\ell(S)=O(\sq), R(S)=O(\sq)$, how many steps
(of the deterministic game) are necessary to reach $\T(\eps(N),1,N)$
where $\eps(N)\to 0$ as $N\to\infty$. From the proof of
Theorem~\ref{det_BS} it can be deduced that if $\eps(N)\sim N^{-\alpha}$,
$0<\alpha<1/2$, then $O(N^{\frac{1}{2}+36\alpha})$ moves suffice.
However, this result is only nontrivial when $\alpha<1/72$
(since $O(N)$ moves are always enough to reach the ``exact''
triangle), and even then it is almost certainly far from
being precise.

Also, loosely speaking, Theorem~\ref{ran_BS} shows that
the typical deviation from the triangle is of order
at most $O(N^{\frac{1}{2}-\frac{1}{144}})$. Again, we do not
believe that that result is the best possible one.
In fact, the author has strong reasons to conjecture
that the typical deviation should be of order $N^{\frac{1}{4}}$;
however, the proof of that is still beyond our reach.

\section*{Acknowledgements}
The author is thankful to Pablo Ferrari for many useful discussions
about the random Bulgarian solitaire,
and to Ira Gessel, who posed the problem of finding the limiting shape
for that model during the Open Problems session at the
conference {\it Discrete Random Walks 2003\/} (IHP, Paris).
Also, the author thanks the anonymous referees for careful reading
of the manuscript and useful comments and suggestions.

\end{document}

%% file: p1.pstex_t
\begin{picture}(0,0)%
\includegraphics{p1.pstex}%
\end{picture}%
\setlength{\unitlength}{4144sp}%
\begingroup\makeatletter\ifx\SetFigFont\undefined
\def\x#1#2#3#4#5#6#7\relax{\def\x{#1#2#3#4#5#6}}%
\expandafter\x\fmtname xxxxxx\relax \def\y{splain}%
\ifx\x\y   
\gdef\SetFigFont#1#2#3{%
  \ifnum #1<17\tiny\else \ifnum #1<20\small\else
  \ifnum #1<24\normalsize\else \ifnum #1<29\large\else
  \ifnum #1<34\Large\else \ifnum #1<41\LARGE\else
     \huge\fi\fi\fi\fi\fi\fi
  \csname #3\endcsname}%
\else
\gdef\SetFigFont#1#2#3{\begingroup
  \count@#1\relax \ifnum 25<\count@\count@25\fi
  \def\x{\endgroup\@setsize\SetFigFont{#2pt}}%
  \expandafter\x
    \csname \romannumeral\the\count@ pt\expandafter\endcsname
    \csname @\romannumeral\the\count@ pt\endcsname
  \csname #3\endcsname}%
\fi
\fi\endgroup
\begin{picture}(4107,3354)(259,-2773)
\put(4366, 74){\makebox(0,0)[lb]{\smash{\SetFigFont{12}{14.4}{rm}{$(7,1)$}%
}}}
\put(4276,-2401){\makebox(0,0)[lb]{\smash{\SetFigFont{12}{14.4}{rm}{$(1,1)$}%
}}}
\put(271,-961){\makebox(0,0)[lb]{\smash{\SetFigFont{12}{14.4}{rm}{$(7,7)$}%
}}}
\end{picture}

%% file: p2.pstex_t
\begin{picture}(0,0)%
\includegraphics{p2.pstex}%
\end{picture}%
\setlength{\unitlength}{4144sp}%
\begingroup\makeatletter\ifx\SetFigFont\undefined
\def\x#1#2#3#4#5#6#7\relax{\def\x{#1#2#3#4#5#6}}%
\expandafter\x\fmtname xxxxxx\relax \def\y{splain}%
\ifx\x\y   
\gdef\SetFigFont#1#2#3{%
  \ifnum #1<17\tiny\else \ifnum #1<20\small\else
  \ifnum #1<24\normalsize\else \ifnum #1<29\large\else
  \ifnum #1<34\Large\else \ifnum #1<41\LARGE\else
     \huge\fi\fi\fi\fi\fi\fi
  \csname #3\endcsname}%
\else
\gdef\SetFigFont#1#2#3{\begingroup
  \count@#1\relax \ifnum 25<\count@\count@25\fi
  \def\x{\endgroup\@setsize\SetFigFont{#2pt}}%
  \expandafter\x
    \csname \romannumeral\the\count@ pt\expandafter\endcsname
    \csname @\romannumeral\the\count@ pt\endcsname
  \csname #3\endcsname}%
\fi
\fi\endgroup
\begin{picture}(5154,7764)(709,-7093)
\put(2206,-1051){\makebox(0,0)[lb]{\smash{\SetFigFont{12}{14.4}{rm}{$(7,2,1,1)$}%
}}}
\put(1036,-3931){\makebox(0,0)[lb]{\smash{\SetFigFont{12}{14.4}{rm}{$(4,6,1)$}%
}}}
\put(4276,-3886){\makebox(0,0)[lb]{\smash{\SetFigFont{12}{14.4}{rm}{$(6,4,1)$}%
}}}
\put(1081,-6721){\makebox(0,0)[lb]{\smash{\SetFigFont{12}{14.4}{rm}{$(3,5,3)$}%
}}}
\put(4321,-6766){\makebox(0,0)[lb]{\smash{\SetFigFont{12}{14.4}{rm}{$(5,3,3)$}%
}}}
\end{picture}

%% file: p3.pstex_t
\begin{picture}(0,0)%
\includegraphics{p3.pstex}%
\end{picture}%
\setlength{\unitlength}{4144sp}%
\begingroup\makeatletter\ifx\SetFigFont\undefined
\def\x#1#2#3#4#5#6#7\relax{\def\x{#1#2#3#4#5#6}}%
\expandafter\x\fmtname xxxxxx\relax \def\y{splain}%
\ifx\x\y   
\gdef\SetFigFont#1#2#3{%
  \ifnum #1<17\tiny\else \ifnum #1<20\small\else
  \ifnum #1<24\normalsize\else \ifnum #1<29\large\else
  \ifnum #1<34\Large\else \ifnum #1<41\LARGE\else
     \huge\fi\fi\fi\fi\fi\fi
  \csname #3\endcsname}%
\else
\gdef\SetFigFont#1#2#3{\begingroup
  \count@#1\relax \ifnum 25<\count@\count@25\fi
  \def\x{\endgroup\@setsize\SetFigFont{#2pt}}%
  \expandafter\x
    \csname \romannumeral\the\count@ pt\expandafter\endcsname
    \csname @\romannumeral\the\count@ pt\endcsname
  \csname #3\endcsname}%
\fi
\fi\endgroup
\begin{picture}(4152,3984)(79,-3223)
\put(4231,-61){\makebox(0,0)[lb]{\smash{\SetFigFont{12}{14.4}{rm}{$\theta_N$}%
}}}
\put(2746,524){\makebox(0,0)[lb]{\smash{\SetFigFont{12}{14.4}{rm}{$V_+$}%
}}}
\put(1561,164){\makebox(0,0)[lb]{\smash{\SetFigFont{12}{14.4}{rm}{$h_+$}%
}}}
\put(796,-1156){\makebox(0,0)[lb]{\smash{\SetFigFont{12}{14.4}{rm}{$V_-$}%
}}}
\put(3091,-376){\makebox(0,0)[lb]{\smash{\SetFigFont{12}{14.4}{rm}{$h_-$}%
}}}
\end{picture}

%% file: p4.pstex_t
\begin{picture}(0,0)%
\includegraphics{p4.pstex}%
\end{picture}%
\setlength{\unitlength}{4144sp}%
\begingroup\makeatletter\ifx\SetFigFont\undefined
\def\x#1#2#3#4#5#6#7\relax{\def\x{#1#2#3#4#5#6}}%
\expandafter\x\fmtname xxxxxx\relax \def\y{splain}%
\ifx\x\y   
\gdef\SetFigFont#1#2#3{%
  \ifnum #1<17\tiny\else \ifnum #1<20\small\else
  \ifnum #1<24\normalsize\else \ifnum #1<29\large\else
  \ifnum #1<34\Large\else \ifnum #1<41\LARGE\else
     \huge\fi\fi\fi\fi\fi\fi
  \csname #3\endcsname}%
\else
\gdef\SetFigFont#1#2#3{\begingroup
  \count@#1\relax \ifnum 25<\count@\count@25\fi
  \def\x{\endgroup\@setsize\SetFigFont{#2pt}}%
  \expandafter\x
    \csname \romannumeral\the\count@ pt\expandafter\endcsname
    \csname @\romannumeral\the\count@ pt\endcsname
  \csname #3\endcsname}%
\fi
\fi\endgroup
\begin{picture}(6054,2274)(79,-1513)
\put(3481,-661){\makebox(0,0)[lb]{\smash{\SetFigFont{12}{14.4}{rm}{$\hh$}%
}}}
\put(4621,-481){\makebox(0,0)[lb]{\smash{\SetFigFont{12}{14.4}{rm}{$m_2$}%
}}}
\put(1996,449){\makebox(0,0)[lb]{\smash{\SetFigFont{12}{14.4}{rm}{$\ha$}%
}}}
\put(5281,-1141){\makebox(0,0)[lb]{\smash{\SetFigFont{12}{14.4}{rm}{$U_2$}%
}}}
\put(4456,494){\makebox(0,0)[lb]{\smash{\SetFigFont{12}{14.4}{rm}{$\ha$}%
}}}
\put(3301,299){\makebox(0,0)[lb]{\smash{\SetFigFont{12}{14.4}{rm}{$U_1$}%
}}}
\put(916,-466){\makebox(0,0)[lb]{\smash{\SetFigFont{12}{14.4}{rm}{$\theta_N$}%
}}}
\put(4111,-1006){\makebox(0,0)[lb]{\smash{\SetFigFont{12}{14.4}{rm}{$\displaystyle\frac{\hh}{5}$}%
}}}
\put(2611,104){\makebox(0,0)[lb]{\smash{\SetFigFont{12}{14.4}{rm}{$\displaystyle\frac{\hh}{5}$}%
}}}
\put(5731,-61){\makebox(0,0)[lb]{\smash{\SetFigFont{12}{14.4}{rm}{$\hh$}%
}}}
\put(3181,-256){\makebox(0,0)[lb]{\smash{\SetFigFont{12}{14.4}{rm}{$m_1$}%
}}}
\end{picture}